\documentclass[fleqn]{mat01}
\usepackage{times,mathtimy,amssymb,latexsym,amsbsy}
\begin{document}

\setcounter{page}{409}
\firstpage{409}

\def\a{\alpha}
\def\b{\beta}
\def\d{\delta}
\def\v{\varepsilon}
\def\r{\rho}
\def\s{\sigma}
\def\t{\tau}
\def\o{\omega}
\def\ov{\overline}
\def\p{\partial}
\def\th{\theta}
\def\ld{\lambda}
\def\pt{P\setminus\{\theta\}}
\def\ph{\varphi}\def\pp{{\prime\prime}}
\def\G{\Gamma}
\def\D{\Delta}
\def\O{\Omega}
\def\e{\eqno}
\def\es{\emptyset}
\def\f{\frac}
\def\ra{\rightarrow}
\def\v{\varepsilon}
\def\l{\leq}
\def\g{\geq}
\def\i{\infty}
\def\su{\sup\limits_{t\in J}}
\def\sq{\sqrt}
\def\st{\stackrel}
\def\fa{\forall}

\def\lms{\lim\limits_}

\def\iii{\int_0^1}
\def\iit{\int_0^t}

\newtheorem{theoree}{Theorem}
\renewcommand\thetheoree{\arabic{section}.\arabic{theoree}}
\newtheorem{theor}[theoree]{\bf Theorem}
\newtheorem{rem}[theoree]{Remark}
\newtheorem{propo}[theoree]{\rm PROPOSITION}
\newtheorem{lem}[theoree]{Lemma}
\newtheorem{definit}[theoree]{\rm DEFINITION}
\newtheorem{coro}[theoree]{\rm COROLLARY}
\newtheorem{exampl}{Example}
\newtheorem{case}{Case}
\newtheorem{pot}[theoree]{Proof of Theorem}

\newtheorem{theore}{\bf Theorem}
\renewcommand\thetheore{\Alph{theore}}

\def\rema{\trivlist \item[\hskip \labelsep{\it Remark.}]}
\def\corol{\trivlist \item[\hskip \labelsep{COROLLARY}]}

\title{Multiple positive solutions to third-order three-point singular
semipositone boundary value problem}

\markboth{Huimin Yu, L Haiyan and Yansheng Liu}{Singular semipositone BVP}

\author{HUIMIN YU, L HAIYAN and YANSHENG LIU$^{*}$}

\address{Department of Mathematics, Shandong Normal University, Jinan 250~014, People's Republic of China\\
\noindent $^{*}$Corresponding Author.\\
\noindent E-mail: ysliu6668@sohu.com}

\volume{114}

\mon{November}

\parts{4}

\Date{MS received 11 April 2004}

\begin{abstract}
By using a specially constructed cone and the fixed point index theory,
this paper investigates the existence of multiple positive solutions for
the third-order three-point singular semipositone BVP:
\begin{equation*}
\begin{cases}
x'''(t)-\ld f(t,x) =0, &t\in(0,1);\\[.3pc]
x(0)=x'(\eta)=x''(1)=0, &
\end{cases}
\end{equation*}
where $\frac{1}{2}<\eta<1$, the non-linear term
$f(t,x)\!\!:(0,1)\times(0,+\i)\to (-\i,+\i)$ is continuous and may be
singular at $t=0$, $t=1$, and $x=0$, also may be negative for some
values of $t$ and $x$, $\ld$ is a positive parameter.
\end{abstract}

\keyword{Singular semipositone boundary value problem; cone; positive
solution; fixed point theorem.}

\maketitle

\section{Introduction}

Singular boundary value problem (BVP) and semipositone BVP arise in a
variety of differential applied mathematics and physics and hence, they
have received much attention (for example, see \cite{1,2,4,6,7} and
references therein). Meanwhile a lot of attention has also been paid to
third-order three-point BVP:
\begin{equation*}
\begin{cases}
 x'''(t)-\lambda f(t,x)=0, &t\in(0,1);\\[.3pc]
 x(0)=x'(\eta)=x''(1)=0, &
\end{cases}\tag{$1_{\lambda}$}
\end{equation*}
where $\frac{1}{2}<\eta<1$, $f(t,x)\!\!:(0,1)\times[0,+\i)\rightarrow(-\i,+\i)$.

In 1998, Anderson \cite{3} considered the problem $(1_\lambda)$ and obtained
an existence result about positive solutions when $f(t,l)=f(l)$ and
$f\hbox{:}\ [0,+\i)\rightarrow[0,+\i)$. Recently, Yao \cite{8} has investigated
$(1_{\lambda})$ when $f$ is semipositone but not singular at $t=0, t=1$,
and $x=0$; and he obtained the following existence theorem.

\begin{theore}[\!]
Suppose
\begin{enumerate}
\renewcommand\labelenumi{{\rm (\arabic{enumi})}}
\leftskip .15pc
\item $\inf\{f(t,l)\!\!:(t,l)\in[0,1]\times [0,+\i)\}=-M>-\i,$ where $M \g 0$.\vspace{.2pc}
\item $B=\max\{f(t,l)\!\!:(t,l)\in[0,1]\times[0,1]\}+M\g 0$.\vspace{.2pc}
\item There exist $0<\a<\b<1$ such that ${\lim}_{l\to +\i}{\min}_{\a\l
t\l\b}\frac{f(t,l)}{l}=+\i.$\vspace{-.5pc}
\end{enumerate}
Then the problem $(1_{\ld})$ has at least one positive solution{\rm ,}
provided
\begin{equation*}
\hskip -1.25pc 0<\lambda<\min\left\{\frac{6}{B\eta^2(3-2\eta)},\quad
\frac{6(2\eta-1)}{M[1-3(1-\eta)^{2}]},\quad \frac{1}{M}\right\}.
\end{equation*}
\end{theore}

By using the approximation method{\rm ,} the fixed point index
theory{\rm ,} and a newly-constructed cone{\rm ,} the present paper
considers $(1_\ld)$ when the non-linear term $f(t,x)$ may be singular at
$t=0,t=1${\rm ,} and $x=0${\rm ,} also may be negative for some values of $t$ and
$x$. The existence of multiple positive solutions is obtained under a
simple assumption which is very similar to that of {\rm \cite{8}}.

The paper is organized as follows{\rm :} In the rest of this
section{\rm ,} some preliminaries are introduced. In \S{\rm 2,} the
main result will be stated and proved and in \S{\rm 3,} some examples
are worked out to demonstrate our main result.\vspace{.5pc}

A map $y\in C[[0,1],R]$ is said to be a positive solution to
BVP$(1_\ld)$\ if it satisfies $(1_\ld)$ and $y(t)>0$ for $t\in (0,1]$.

In obtaining positive solutions to $(1_\ld)$, the following two results
are fundamental.

\begin{lem}\hskip -.3pc {\rm (\cite{5},} Lemma {\rm 2.3.1,} p.~{\rm 88,} and Lemma
{\rm 2.3.2,} p.~{\rm 91)}.\ \ Let $P$ be a cone of real Banach space $E, \O$
be a bounded open set of $E, \theta\in \O${\rm ,} the map
$A\!\!: P\cap\overline{\O}\rightarrow P$ be completely continuous.

\begin{enumerate}
\renewcommand\labelenumi{{\rm \arabic{enumi})}}
\leftskip -.15pc
\item If $x\neq \mu Ax$ for $x\in P\cap\partial\O$ and
$\mu\in[0,1]${\rm ,} then $i(A,P\cap\O,P)=1${\rm ;}\vspace{.2pc}
\item If ${\inf}_{x\in P\cap\partial\O}\|Ax\|>0$ and $Ax\neq\mu x$ for
$x\in P\cap\partial\O$ and $\mu \in (0,1]${\rm ,} then
$i(A,P\cap\O,P)=0$.
\end{enumerate}
\end{lem}

\begin{lem}\hskip -.3pc {\rm (\cite{6},} Lemma {\rm 1.2)}.\ \
If $g\in C[(0,+\i),R^{+}]${\rm ,} then there exists a non-decreasing
function $h\in C[R^{+},R^{+}]$ such that $h(x)>0$ as $x>0$ and
$g(x)h(x)\in C[R^{+},R^{+}]$ {\rm (}that is{\rm ,} ${\lim}_{x\to
0^+}g(x)h(x)$ exists{\rm ),} where $R^{+} = [0,+\i)$.
\end{lem}

\section{Main results}

For convenience, we list the following hypothesis.

\noindent $(\hbox{H}_{1})$ There exists $M>0$ such that
\begin{equation*}
0\l M + f(t,x)\l g(t)h(x),\quad \fa t\in (0,1),\quad x\in(0,+\i),
\end{equation*}
where $g(t)\!\!:(0,1)\rightarrow(0,+\i), \int_0^{\eta} sg(s)\hbox{d} s
+ \int_{\eta}^{1}g(s)\hbox{d}s<+\i,$ and $h\!\!:(0,+\i)\rightarrow [0,+\i)$
is continuous.

\noindent $(\hbox{H}_{2})$ There exists $[\a,\b]\subset(0,1)$ such that
\begin{equation*}
\lim_{x\rightarrow 0+}f(t,x) = +\i ,\quad  \lim_{x\to +\i}\frac{f(t,x)}{x}=+\i,
\end{equation*}
both uniformly with respect to $t\in [\a,\b]$.

The following theorem is our main result.
\setcounter{theoree}{0}
\begin{theor}[\!]
Suppose that conditions $(\hbox{H}_{1})$ and $(\hbox{H}_{2})$ hold{\rm ,} then for each
$r>0${\rm ,} there exists $\overline{\ld}=\overline{\ld}(r)>0$ such that BVP
$(1_\ld)$ has at least two positive solutions $x(t)$ and $y(t)$ in
$C^{3}[(0,1),R]\cap C[[0,1],R]$ satisfying $0<\|x\|<r<\|y\|$ provided
$\ld \in (0,\overline{\ld})$.
\end{theor}

Before giving the proof of Theorem~2.1, we first list some preliminaries
and prove some lemmas.

Let $I = [0,1], E = C[I,R]$, then $E$ is a Banach space with norm
$\|x\|={\max}_{t\in I}|x(t)|$. Throughout this paper, we shall
use the following notation:
\begin{equation*}
G(t,s)=\begin{cases}
ts-\frac{1}{2}t^{2},    &0\l s\leq \eta,0\l t\l s;\\[.3pc]
\frac{1}{2}s^{2},   &0\l s\leq \eta,0\l s\l t;\\[.3pc]
\eta t-\frac{1}{2}t^2,  &\eta \l s\l1,0\l t\l s;\\[.3pc]
\frac{1}{2}s^{2}-ts+\eta t, &\eta \l s\l1,0\l s\l t.
\end{cases}
\end{equation*}

It is well-known that $G(t,s)$ is the Green's function of
homogeneous boundary value problem:
\begin{equation*}
\begin{cases}
x'''(t)=0, &0\l t\l 1;\\[.3pc]
x(0)=x'(\eta)=x''(1)=0.\end{cases}
\end{equation*}

\setcounter{theoree}{0}
\begin{lem}
$G(t,s)$ defined as above have the following properties{\rm :}
\begin{enumerate}
\renewcommand\labelenumi{{\rm (\arabic{enumi})}}
\leftskip .15pc
\item $J(s) =: \mathop{\max}\limits_{t\in I}G(t,s)=\begin{cases}
\frac{1}{2}s^2, &0\l s\l\eta;\\[.3pc]
\frac{1}{2}\eta^2,\quad \eta \l s\l 1.
\end{cases}$\vspace{.2pc}

\item $G(t,s)\g q(t)J(s)${\rm ,}\ \,where\ \,$q(t) = \begin{cases}
\eta t, &0\l t\l\eta;\\[.3pc]
2\eta t-t^2, &\eta\l t\l1,\end{cases}$ is a non-negative concave
function in I.\vspace{.2pc}

\item $\Psi^{*}(t)=:\iii G(t,s){\rm d}s = \frac{1}{6}[t^3-3t^2+(6\eta
-3\eta^2)t]$ and $\|\Psi^{*}\|= {\max}_{t\in I}\Psi^{*}(t)=\frac{1}{6}(3\eta^{2} - 2\eta^{3}).$\vspace{.2pc}

\item $\Psi^{*}(t)\l Kq(t), \fa t\in I,$ where $K = \max\{1,B_0\}, B_0=(6\eta-3\eta^2-2)[6(2\eta-1)]^{-1}$.
\end{enumerate}
\end{lem}

\begin{proof}
The proof of (1) and (2) can be seen from \cite{8}.

The proof of (3) is as follows: For $t\in [\eta,1]$, we can get
\begin{align*}
\iii G(t,s)\hbox{d}s &= \int_{0}^{\eta} \frac{1}{2}s^{2}\hbox{d}s + \int_\eta^{t}\left(\frac{1}{2}s^{2} - ts + \eta t\right) \hbox{d}s + \int_{t}^{1}\left(\eta t - \frac{1}{2}t^{2}\right)\hbox{d}s\\[.3pc]
&= \frac{1}{6}[t^3-3t^2+(6\eta -3\eta^2)t].
\end{align*}
Meanwhile, for $t\in [0,\eta]$,
\begin{align*}
\iii G(t,s)\hbox{d}s &= \int_{0}^{t} \frac{1}{2}s^{2}\hbox{d}s + \int_t^\eta\left(ts - \frac{1}{2}t^2\right)\hbox{d}s + \int_\eta^1\left(\eta t - \frac{1}{2}t^{2}\right)\hbox{d}s\\[.3pc]
&= \frac{1}{6}[t^3-3t^2+(6\eta -3\eta^2)t].
\end{align*}

Therefore it is easy to see that
\begin{equation*}
\mathop{\max}\limits_{t\in I} \Psi^*(t) = \Psi^*(\eta) = \frac{1}{6}(3\eta^2-2\eta^3).
\end{equation*}

The proof of (4) is as follows: For $t\in [0,\eta]$, we have
\begin{equation*}
\Psi^*(t)-Kq(t)\l\Psi^*(t)-q(t) = \frac{1}{6}(t^3-3t^2-3\eta^2t).
\end{equation*}

Let $p(t) = \frac{1}{6}(t^2-3t-3\eta^2)$. Noticing that $p(0)<0,
p(\eta)<0$, we assert that $p(t)<0$ for $ t\in [0,\eta]$. And follows
the conclusion $\Psi^*(t)-Kq(t)\l0$ for $t\in [0,\eta]$.

On the other hand, for $t\in [\eta, 1]$,
\begin{align*}
\Psi^*(t)-Kq(t) &\l \Psi^*(t)-B_0q(t)\\[.3pc]
 &= \frac{1}{6}[t^{3}-3t^{2} + (6\eta - 3\eta^{2})t]- \frac{6\eta-3\eta^2-2}{6(2\eta-1)}(2\eta t-t^2).
\end{align*}

Similar to the above, for $t\in [\eta,1]$, let
\begin{equation*}
p(t)= \frac{1}{6}(t^2-3t+6\eta -3\eta^2) -
\frac{6\eta-3\eta^2-2}{6(2\eta-1)}(2\eta -t).
\end{equation*}

Since $p(\eta)=(-\eta^3-\eta+2\eta^2)[6(2\eta -1)]^{-1} <0$\ and \
$p(1)=0$, we can prove that $\Psi^*(t)-Kq(t)\l 0$ for $t\in [\eta,1].$
\hfill $\Box$
\end{proof}

Let $P=\{x\in E\!\!:\ x(t)\g 0,\ \fa t\in I\},\ Q=\{x\in P\!\!:\ x(t)\g
q(t)\|x\|,\ \fa t\in I\}$. Obviously $P$ and $Q$ are cones in $E$.

For $\ld\in (0,+\i), j\in N\ (N=\{1,2,3,\ldots\})$,\ consider the
following approximation problem of $(1_\ld)$:
\setcounter{equation}{1}
\begin{equation}
\begin{cases}
x'''(t)-\ld f_{j}^{*} (t,x(t)-\phi_{\ld}(t)+\frac{1}{j}) = 0, &t\in(0,1);\\[.3pc]
x(0)=x'(\eta)=x''(1)=0, &
\end{cases}
\end{equation}
where
\begin{align*}
\phi_{\ld}(t) &= \ld M\Psi^*(t),\\[.3pc]
f_{j}^{*}\left(t,u+\frac{1}{j}\right) &= \begin{cases}
f (t,u + \frac{1}{j}) + M, &u\g0;\\[.3pc]
f (t,\frac{1}{j}) + M, &u<0 .
\end{cases}
\end{align*}

It is easy to see $\|\phi_{\ld}\|=\ld
M\|\Psi^*\| = \frac{\ld}{6}M(3\eta^{2}-2\eta^{3})<+\,\i$. Evidently, $x\in
C^3[(0,1),R]\cap C[I,R]$\ is a solution of (2) if and only if $x\in
C[I,R]$ is a solution of the following integral equation:
\begin{equation}
x(t)= \ld\iii G(t,s)f_j^{*}\left(s,x(s)-\phi_\ld(s)+\frac{1}{j}\right)\hbox{d}s.
\end{equation}
Let operator $A_\ld^{j}$ be defined by
\begin{equation}
(A_\ld^{j} x)(t)=:\ld\iii G(t,s)f_j^{*}\left(s,x(s)-\phi_\ld(s)+\frac{1}{j}\right)\hbox{d}s,\quad \fa t\in I.
\end{equation}

\begin{lem}
For each $\ld \in (0,+\i), j\in N, A_\ld^j\!\!:Q\rightarrow Q$ is a continuous
and compact mapping.
\end{lem}

\begin{proof}
From Lemma~2.1(2), we can see that $A_\ld^j\!\!:Q\rightarrow Q$. Meanwhile,
by $f_j^*\in C[(0,1)\times R,R^+]$ and $(\hbox{H}_1)$, one can conclude that
$A_\ld^j$ is continuous and compact from $Q$\ to $Q$.\hfill$\Box$
\end{proof}

Therefore, $x\in Q$ is a solution of (2) if $x$ is a fixed point of
$A_\ld^j$ on $Q$. Then we next consider the existence of fixed point of
$A_\ld^j$ on $Q$.

\begin{lem}
For each $r>0${\rm ,} there exists $\ld(r)>0$
such that
\begin{equation*}
i(A_\ld^j,Q_r,Q)=1,\qquad \fa \ld \in (0,\ld(r)), j\in N,
\end{equation*}
where $Q_r=\{x\in Q:\|x\|<r\}$.
\end{lem}

\begin{proof}
By Lemma 1.2 and $(\hbox{H}_{1})$, there exists a non-decreasing function $p\in
C[R^+,R^+]$ such that $p(x)h(x)\in C[R^+,R^+]$ and $p(x)>0$ for
$\fa x > 0$. For $\fa r>0,$ let
\begin{equation*}
\ld (r)=:\min\left\{\frac{r}{2KM},\frac{b}{c(r)a}\right\},
\end{equation*}
where
\begin{align*}
a &=: \frac{1}{2}\int_0^\eta \t^2g(\t)\hbox{d}\t + \frac{\eta^{2}}{2}\int_\eta^1
g(\t) \hbox{d}\t,\quad c(r)=:\mathop{\max}\limits_{x\in [0,r+1]}h(x)p(x),\\[.3pc]
b &=: 2\int_0^{\frac{r(2\eta-1)}{2}}p(s)\hbox{d}s+2rp\left(\frac{r(2\eta-1)}{2}\right)
(1-\eta),\\[.3pc]
&\quad\,K\ \hbox{is defined in Lemma 2.1(4).}
\end{align*}

We now claim that
\begin{equation*}
x\neq\mu(A_\ld^j x),\quad
\fa \mu\in [0,1],\quad x\in \partial Q_r,\quad \ld \in (0,\ld(r)).
\end{equation*}

If it is false, then there exist $x_0\in \partial Q_r$ and
$\mu_0\in[0,1]$ such that
\begin{equation*}
x_0(t)=\mu_0(A_\ld^j x_0(t)),\quad \fa t\in I.
\end{equation*}

Noticing $x_0 \in Q$\ and using Lemma 2.1, we get
\begin{align*}
x_0(t)\g q(t)\|x_0\| &= rq(t),\\[.3pc]
\phi_\ld (t) &= \ld M\Psi^*(t)\l\ld MKq(t)\l \frac{\ld KM}{r}x_0(t),\quad \fa t\in I.
\end{align*}

So
\begin{equation}
x_0(t)-\phi_\ld(t)\g \left(1 - \frac{\ld KM}{r}\right)x_0(t)\g
\frac{1}{2}x_0(t)\g \frac{r}{2}q(t),\quad \fa t\in I.
\end{equation}
Furthermore, since $q(t)$ is a non-negative function on [0,1], it
is easy to see
\begin{equation}
f_j^*\left(t,x_0(t)-\phi_\ld(t)+\frac{1}{j}\right) = f\left(t,x_0(t)-\phi_\ld
(t)+\frac{1}{j}\right)+M,\quad \fa t\in I.
\end{equation}
Therefore, by (5), (6), and $(\hbox{H}_1)$ we have
\begin{align}
x_0'''(t) &= \ld \mu_0 \left[f\left(t,x_0(t)-\phi_\ld (t)+\frac{1}{j}\right)+M\right]\nonumber\\[.3pc]
&\l \ld g(t)h\left(x_0(t)-\phi_\ld (t)+\frac{1}{j}\right)\l \ld
g(t)\frac{c(r)}{p(x_0(t)-\phi_\ld(t))}\nonumber\\[.3pc]
&\l \ld g(t)\frac{c(r)}{p (\frac{1}{2}x_0(s))},\quad \fa t\in (0,1).
\end{align}

Integrate (7) from $t$ to 1 to obtain
\begin{align}
-x_0''(t) &\l \ld c(r)\int
_t^1\frac{g(s)}{p (\frac{1}{2}x_0(s))}\,\hbox{d}s\nonumber\\[.3pc]
&\l \ld c(r)\int_t^1\frac{g(s)}{p (\frac{1}{2}\min
\{x_0(t),x_0(1)\})}\,\hbox{d}s,\quad t\in(0,1).
\end{align}
Now integrating (8) from $t$ to $\eta$, we have
\begin{align*}
x_0'(t) &\l \ld c(r)\int_t^\eta\int_s^1\frac{g(\t)}{p (\frac{1}{2}\min
\{x_0(s),x_0(1)\})}\,\hbox{d}\t \hbox{d}s\\[.3pc]
&\l \ld c(r)\int_t^\eta\int_s^1 \frac{g(\t)}{p (\frac{1}{2}\min
\{x_0(t),x_0(1)\})}\,\hbox{d}\t \hbox{d}s,\qquad \fa t\in (0,\eta),
\end{align*}
that is,
\begin{equation}
x_0'(t)p\left(\frac{1}{2}\min
\{x_0(t),x_0(1)\}\right)\l \ld c(r)\int _t^\eta\int_s^1 g(\t)\hbox{d}\t \hbox{d}s,\quad
t\in(0,\eta).
\end{equation}
Integrate (9) from 0 to $\eta$ to get
\begin{align*}
&\int_0^\eta x_0'(t)p\left(\frac{1}{2}\min \{x_0(t),x_0(1)\}\right)\hbox{d}t\l \ld
c(r)\int_0^\eta\int _t^\eta\int_s^1g(\t)\hbox{d}\t \hbox{d}s\hbox{d}t\\[.3pc]
&\quad\,\l \ld c(r)\left[\int_\eta^1g(\t)\hbox{d}\t\int_0^\eta \hbox{d}s\int_0^s \hbox{d}t+\int
_0^\eta g(\t)\hbox{d}\t\int _0^\t \hbox{d}s\int_0^s \hbox{d}t\right]\\[.3pc]
&\quad\,=\ld c(r)\left[\frac{\eta^{2}}{2}\int_\eta^1g(\t)\hbox{d}\t+\frac{1}{2}\int
_0^\eta \t^2g(\t)\hbox{d}\t\right]=\ld c(r)a.
\end{align*}
On the other hand, noticing $x_0\in Q$ and $x_0'(t)>0$ for $t\in
(0,\eta)$, we get
\begin{align*}
&\int _0^\eta x_0'(t)p\left(\frac{1}{2}\min \{x_0(t),x_0(1)\}\right)\hbox{d}t\\[.3pc]
&\quad\,\g\int_0^\eta x_0'(t)p\left(\frac{1}{2}\min \{x_0(t),r(2\eta-1)\}\right)\hbox{d}t\\[.3pc]
&\quad\,=\int_0^r p\left(\frac{1}{2}\min \{s,r(2\eta-1)\}\right)\hbox{d}s.
\end{align*}
Then
\begin{align*}
&\int _0^r p\left(\frac{1}{2}\min \{s,r(2\eta-1)\}\right)\hbox{d}s\\[.3pc]
&\quad\,=\int_0^{r(2\eta-1)}p\left(\frac{1}{2}s\right)\hbox{d}s+\int_{r(2\eta-1)}^{r}p\left(\frac{r(2\eta-1)}{2}\right)\hbox{d}s\\[.3pc]
&\quad\,=2\int_0^{\frac{r(2\eta-1)}{2}}p(s)\hbox{d}s+2rp\left(\frac{r(2\eta-1)}{2}\right)(1-\eta)=b\l\ld c(r)a.
\end{align*}
This guarantees that $\ld\g\frac{b}{c(r)a}$, which is in contradiction
with $0<\ld<\frac{b}{c(r)a}$. Consequently the result of Lemma~2.3
follows.\hfill$\Box$
\end{proof}

\begin{lem}
For $\fa \ld\in (0,\ld(r))${\rm ,} there exists $R > r > 0$ such that
\begin{equation*}
i(A_\ld^j,Q_R,Q)=0,\quad \fa j\in N.
\end{equation*}
\end{lem}

\begin{proof}
By $(\hbox{H}_{2})$ we know that, for $\fa \ld \in (0,\ld(r)),$ there exists
$L>r>0$ such that
\begin{equation}
\frac{f(t,x)}{x}\g 2\left(\ld \mathop{\min}\limits_{t\in [\a,\b]}q(t)\int_\a^\b G(t,s)\hbox{d}s\right)^{-1},\quad \fa x>L, t\in [\a,\b].
\end{equation}
Set
\begin{equation*}
R = R(\ld)>\max\left\{\frac{2L}{\mathop{\min}\limits_{t\in [\a,\b]}q(t)},r\right\}.
\end{equation*}
We now claim that
\begin{equation*}
A_{\ld}^jx\neq\mu x,\quad \fa x\in \partial Q_R,\quad \mu \in (0,1].
\end{equation*}
In fact, if it is not true, then there exist $x_0\in \partial Q_R$ and
$\mu_0 \in (0,1]$ such that $A_\ld^jx_0=\mu_0x_0,$ that is,
\begin{equation}
x_0(t)\g(A_{\ld}^jx_0)(t)=\ld \iii
G(t,s)f_j^{*}\left(s,x_0(s)-\phi_\ld(s)+\frac{1}{j}\right)\hbox{d}s,
\quad \fa t\in I.
\end{equation}
Since
\begin{equation*}
x_0(t)\g q(t)R,\quad \phi_\ld (t)=\ld M\Psi^*(t)\l\ld
MKq(t),\quad \fa t\in I,
\end{equation*}
then for $\fa t\in [\a,\b]$, we have
\begin{equation*}
x_0(t)-\phi_\ld(t)\g(R-\ld MK)\mathop{\min}\limits_{t\in [\a,\b]}
q(t)\g \frac{R}{2}\mathop{\min}\limits_{t\in [\a,\b ]}q(t)>L.
\end{equation*}
By (10) and (11) we get for $\fa t\in [\a,\b]$,
\begin{align*}
x_0(t) &\g \ld \iii G(t,s)f_j^{*}\left(s,x_0(s)-\phi_\ld(s) + \frac{1}{j}\right)\hbox{d}s\\[.3pc]
&=\ld \iii G(t,s)\left[f\left(s,x_0(s)-\phi_\ld (s) + \frac{1}{j}\right)+M\right]\hbox{d}s\\[.3pc]
&>\ld \int_\a^\b G(t,s)f\left(s,x_0(s)-\phi_\ld (s) + \frac{1}{j}\right)\hbox{d}s\\[.3pc]
&\g 2 \ld\left(\ld\mathop{\min}\limits_{t\in [\a,\b]}q(t)\int_\a^\b
G(t,s)\hbox{d}s\right)^{-1}(R-\ld MK)\\[.3pc]
&\quad\,\times \mathop{\min}\limits_{t\in [\a,\b ]}q(t)\int_\a^\b G(t,s)\hbox{d}s\\[.3pc]
&\g 2(R-\ld KM)>R.
\end{align*}
This is a contradiction with $x_0\in \partial Q_R$. Thus, the proof of
Lemma 2.4 is completed.\hfill $\Box$
\end{proof}

\begin{lem}
For the above-mentioned $r>0$ and $\ld(r)>0${\rm ,} there exists
$\overline{\ld}=\overline{\ld}(r)\in (0,\ld(r)]$ satisfying that for
each $\ld \in (0,\overline{\ld})${\rm ,} there exists $r'=r'(\ld)\in (0,r)$
such that $i(A_\ld^j,Q_{r'},Q)=0${\rm ,} if $j$ is sufficiently large.
\end{lem}

\begin{proof}
By condition $(\hbox{H}_{2})$, we can get that for each
$L'>MK ({\min}_{t\in [\a,\b ]}\int_\a^\b
G(t,s)\hbox{d}s)^{-1}$, there exists $\delta >0$ such that
\begin{equation*}
f(t,x)>L'\quad\hbox{for}\quad x\in (0,\d)\quad \hbox{and}\quad t\in [\a,\b].
\end{equation*}
Let
\begin{equation*}
\overline{\ld}=\min\left\{\frac{\d}{M},\ld (r)\right\},\quad
l=L'\mathop{\min}\limits_{t\in [\a,\b ]}\int_\a^\b G(t,s)\hbox{d}s.
\end{equation*}
Then $M\ld<\d$ for $\ld \in (0,\overline{\ld})$. Fix $\ld \in
(0,\overline{\ld})$, choose $r'=r'(\ld)\in (0,\d)$ and $j$ sufficiently
large such that
\begin{equation*}
\frac{r'}{l} < \ld < \frac{r'}{MK},\quad r' + \frac{1}{j} < \d.
\end{equation*}
Next, we prove
\begin{equation*}
A_\ld^jx\neq\mu x,\quad \fa x \in\partial Q_{r'},\quad \mu \in (0,1].
\end{equation*}
Suppose this is false, then there exist $x_0\in \partial Q_{r'}$ and
$\mu_0\in (0,1]$ such that $A_\ld^jx_0=\mu_0x_0$. Since $x_0(t)\g
r'q(t)$ and $\phi_\ld\l\ld MKq(t)$ for $t\in I${\rm ,} we can get
\begin{equation*}
x_0(t)-\phi_\ld(t)\g(r'-\ld MK)q(t)\g 0,\quad \fa t\in I.
\end{equation*}
Consequently
\begin{equation*}
f_j^{*}\left(t,x_0(t)-\phi_\ld(t) + \frac{1}{j}\right) = f\left(t,
x_0(t)-\phi_\ld (t) + \frac{1}{j}\right) + M,\quad \fa t\in I.
\end{equation*}
On the other hand,
\begin{equation*}
0 < x_0(t)-\phi_\ld(t) + \frac{1}{j}\l x_0(t) + \frac{1}{j}\l r' +
\frac{1}{j}<\d,\quad \fa t\in I.
\end{equation*}
Therefore,
\begin{align*}
x_0(t)\g(A_\ld^jx_0)(t) &= \ld \iii
G(t,s)f_j^{*}\left(s,x_0(s)-\phi_\ld(s) + \frac{1}{j}\right)\hbox{d}s\\[.3pc]
&>\ld \iii G(t,s)f\left(s,x_0(s)-\phi_\ld (s) + \frac{1}{j}\right)\hbox{d}s\\[.3pc]
&>\ld L'\int_\a^\b G(t,s)\hbox{d}s \g\ld L'\mathop{\min}\limits_{t\in [\a,\b
]}\int_\a^\b G(t,s)\hbox{d}s>r',\\[.3pc]
&\quad\, \fa t\in [\a,\b].
\end{align*}
This is in contradiction with $x_0\in \partial Q_{r'}$ and immediately our
result follows.\hfill $\Box$
\end{proof}

\begin{lem}
For each $\ld \in (0,\overline{\ld})$ and sufficiently large $j${\rm ,}
BVP {\rm (}2{\rm )} has at least two positive solutions $x_j$ and $y_j$
satisfying
\begin{equation*}
r'<\|x_j\|<r<\|y_j\|< R,
\end{equation*}
where $R, r'$ and $\overline{\ld}$ are the same as in Lemmas {\rm 2.4} and
{\rm 2.5}.
\end{lem}

\begin{proof}
By Lemmas 2.3--2.5 and the additivity of fixed point index we can easily
get
\begin{align*}
i(A_\ld^j,Q_R\setminus \overline{Q_r},Q) &= i(A_\ld^j,Q_R,Q)-i(A_\ld^j,Q_r,Q)=0-1=-1,\\[.3pc]
i(A_\ld^j,Q_r\setminus \overline{Q_{r'}},Q) &= i(A_\ld^j,Q_r,Q)-i(A_\ld^j,Q_{r'},Q)=1-0=1.
\end{align*}

It follows from solution property of the fixed point index that there
exist $y_j \in Q_R\setminus \overline{Q_r}$ and $x_j\in Q_r\setminus
\overline{Q_{r'}}$ such that
\begin{equation*}
A_\ld^jx_j=x_j\quad\hbox{and}\quad A_\ld^jy_j=y_j.
\end{equation*}
Consequently, Lemma 2.6 follows.\hfill $\Box$
\end{proof}

\setcounter{theoree}{0}
\begin{pot}
{\rm Let $\{x_j\},\{y_j\} \ (j\g j_0)$ be the positive solutions obtained in
Lemma~2.6.

We now prove $\{x_j\}_{j\g j_0}$ is a bounded and equicontinuous family
on [0,1].

The boundedness is obvious. Since $\|x_j\|>r', \phi_\ld(t)\l\ld
MKq(t)\l\frac{\ld MK}{r'}x_j(t)$ for $j\g j_0,$ and $\ld MK<r',$ we
have
\begin{equation*}
x_j(t)-\phi_\ld (t)\g \left(1-\frac{\ld MK}{r'}\right)x_j(t),\quad t\in (0,1).
\end{equation*}
Similar to the proof of (9), we get
\begin{equation}
x_j'(t)p\left(\frac{1}{2}\min \{x_j(t), x_j(1)\}\right)\l \ld c(r)
\int_t^\eta\int_s^1 g(\t)\hbox{d}\t \hbox{d}s,\quad \fa t\in (0,\eta).
\end{equation}
Let $I\!\!:{R}^{+}\rightarrow R^{+}$ be defined by
\begin{equation*}
I(u)=:\int_0^up\left(\frac{1}{2}\min\{s,r'(2\eta-1)\}\right)\hbox{d}s.
\end{equation*}
By $(\hbox{H}_{1})$ we get
\begin{equation*}
\zeta=:\int_0^\eta \t g(\t)\hbox{d}\t <+\i,\quad
\gamma=:\int_\eta^1 g(\t)\hbox{d}\t<+\i,
\end{equation*}
and for $\fa \varepsilon>0 $, there exist $\d\in (0,\varepsilon)$ such
that
\begin{equation*}
\int_{t_1}^{t_2}\t g(\t)\hbox{d}\t <\varepsilon\quad \hbox{as}\quad
0\l t_1\l t_2\l\eta\quad \hbox{and}\quad |t_2-t_1|<\d.
\end{equation*}
For $\fa t_1, t_2\in [0,\eta]$ with $|t_2-t_1|<\d$ and $t_1\l t_2$, we
have
\begin{align*}
&|I(x_j(t_2))-I(x_j(t_1))|\\[.3pc]
&\quad\, = \left|\int_{t_1}^{t_2}p\left(\frac{1}{2}\min
\{x_j(t),r'(2\eta-1)\}\right)x_j'(t)\hbox{d}t\right|\\[.3pc]
&\quad\,\l \int_{t_1}^{t_2}p\left(\frac{1}{2}\min\{x_j(t),x_j(1)\}\right)x_j'(t)\hbox{d}t\\[.3pc]
&\quad\,\l\ld c(r)\int_{t_1}^{t_2}\int_t^\eta\int _s^1 g(\t)\hbox{d}\t \hbox{d}s \hbox{d}t\\[.3pc]
&\quad\,=\ld c(r)\left[\frac{1}{2}\int_{t_1}^{t_2}g(\t)(\t-t_1)^2\hbox{d}\t + \frac{1}{2}
(t_2-t_1)^{2}\int_{t_{2}}^{1} g(\t)\hbox{d}\t\right.\\[.3pc]
&\qquad\,\left.+(t_2-t_1)\int _{t_2}^\eta
(\t-t_2)g(\t)\hbox{d}\t +(t_2-t_1)(\eta -t_2)\int_\eta^1 g(\t)\hbox{d}\t\right]
\end{align*}
\begin{align*}
&\quad\,\l \ld c(r) \left[ \eta\int_{t_1}^{t_2}g(\t)\t \hbox{d}\t + \frac{3}{2}
(t_2-t_1) \int_0^ \eta \t g(\t)\hbox{d}\t \phantom{\int_\eta^1}\right.\\[.3pc]
&\qquad\, \left.+ 2\eta (t_2-t_1)\int_\eta^1
g(\t)\hbox{d}\t\right]\\[.3pc]
&\quad\,\l \ld c(r)\left(\eta + \frac{3}{2}\zeta+2\eta \gamma\right)\varepsilon.
\end{align*}
Thus, we obtain that $\{I(x_j(t))\}_{j\g j_0}$ is equicontinuous on $[0,\eta]$.
Furthermore, by the uniform continuity of $I^{-1}$ on
$[0,I(r)]$ and
\begin{equation*}
|x_j(t_2)-x_j(t_1)|=|I^{-1}(I(x_j(t_2)))-I^{-1}(I(x_j(t_1)))|,
\end{equation*}
the equicontinuity of $\{x_j(t)\}_{j\g j_0}$ on $[0,\eta]$ is established.

For $\fa t_0\in (0,\eta), t\in [t_0,1]$, noticing $x_j \in Q$
and returning to (8) (replacing $x_0$ with $x_j$), we obtain
\begin{align}
0\l-x_j''(t) &\l \ld c(r)\int_t^{1}\frac{g(s)}{p (\frac{1}{2}\min
\{x_j(t),x_j(1)\})}\hbox{d}s\nonumber\\[.3pc]
&\l \ld c(r)\int _{t_0}^1\frac{g(s)}{p (\frac{1}{2}\min
\{x_j(t_0),x_j(1)\})}\hbox{d}s\nonumber\\[.3pc]
&\l \frac{\ld c(r)}{p (\frac{1}{2}\min \{r'\eta t_0,
r'(2\eta-1)\})}\int_{t_0}^1 g(s)\hbox{d}s,\quad t\in[t_0,1].
\end{align}

Since
\begin{equation*}
\int_{t_0}^1 g(s)\hbox{d}s = \int_{t_0}^\eta g(s)\hbox{d}s+\int_\eta^1
g(s)\hbox{d}s\l\frac{1}{{t_0}}\int_0^\eta s g(s)\hbox{d}s + \int_\eta^1
g(s)\hbox{d}s,
\end{equation*}
and by $(\hbox{H}_{1})$, we can get $d=:\int_{t_0}^1 g(s)\hbox{d}s<+\i$.

This together with (13) guarantees that
\begin{equation}
|x_j''(t)| \l\frac{\ld c(r)}{p (\frac{1}{2}\min \{r'\eta t_0, r'(2\eta-1)\})}d<\i,\quad
t\in [t_{0}, 1].
\end{equation}
Therefore,
\begin{align*}
|x_j'(t)| &= |x_j'(t)-x_j'(\eta)|=|x_j''(\xi)(t-\eta)|\\[.3pc]
&\l |x_j''(\xi)|\l\frac{\ld c(r)}{p(\frac{1}{2}\min
\{r'\eta t_0, r'(2\eta-1)\})}d,\quad \fa t\in [t_0,1].
\end{align*}

Consequently, the equicontinuity of $\{x_j\}_{j\g j_0}$ on $[t_0,1]$
follows, and now $\{x_j\}_{j\g j_0}$ is equicontinuous in $I$.

The Arzela--Ascoli theorem guarantees the existence of a
subsequence $N_0$ of $N$ and a function $x^*\in C[0,1]$ with $x_j$
converging uniformly on $I$ to $x^*$ as $j\rightarrow +\i$ through $N_0$.
Also $x^*(0)=0,\ r'\l\|x^*\|\l r,$ and $x^*(t)\g r'q(t)$ for $t\in
I.$ In particular, $x^*(t)>0$ on (0,1).

Noticing (14), $\{x_j''(\eta)\}_{j\in {N_0}}$ is a bounded sequence.
Therefore, $\{x_j''(\eta)\}_{j\in {N_0}}$ has a convergent subsequence.
For convenience, let $\{x_j''(\eta)\}_{j\in {N_0}}$ denote this
subsequence also, and let $r_0\in \ {\bf R}$ be its limit.

Fix $t\in [0,1], x_j\  (j\in N_0)$ satisfies the integral
equation
\begin{align}
x_j(t)&=x_j(\eta)+ \frac{(t-\eta)^{2}}{2}x_j''(\eta)\nonumber\\[.3pc]
&\quad\,+\ld \int _\eta^t\frac{(t-s)^{2}}{2}\left[f\left(s,x_j(s)-\phi_\ld(s) + \frac{1}{j}\right) + M\right]\hbox{d}s.
\end{align}
Now, let $j\to +\i$ through $N_0$ in (15) to obtain
\begin{align}
x^*(t) &= x^*(\eta) + \frac{(t-\eta)^{2}}{2}r_0\nonumber\\[.3pc]
&\quad\,+\ld \int_\eta^{t}\frac{(t-s)^{2}}{2}[f(s,x^*(s)-\phi_\ld(s))+M
]\hbox{d}s,\quad t\in I.
\end{align}
Differentiate (16) directly to obtain ${x^*}'(\eta)=0$. On the other
hand, $x_{j}\ (j\in N_0)$ also satisfies the integral equation
\begin{align}
x_j(t)&=x_j(1)-x_j'(1)(1-t)\nonumber\\
&\quad\,-\ld \int_t^1\frac{(s -t)^{2}}{2}\left[f\left(s,
x_j(s)-\phi_\ld(s)+\frac{1}{j}\right)+M\right]\hbox{d}s,\quad t\in I.
\end{align}
This together with $(\hbox{H}_{1})$ and
$r'\l\|x_j\|\l r$ guarantees that $\{x_j'(1)\}_{j\in {N_0}}$ is a
bounded sequence. Therefore, $\{x_j'(1)\}_{j\in {N_0}}$ has a
convergent subsequence. Also let $\{x_j'(1)\}_{j\in {N_0}}$ denote
this subsequence, and let $R_0\in {\bf R}$ be its limit.

Let $j\rightarrow +\i$ through $N_0$ in (17) to get
\begin{equation}
x^*(t)=x^*(1)-R_0(1-t)-\ld \int_t^1\frac{(s
-t)^{2}}{2}[f(s,x^*(s)-\phi_\ld(s))+M]\hbox{d}s.
\end{equation}
Differentiate (18) directly to obtain ${x^*}''(1)=0$. Therefore,
$x^*(t)$ is a solution of
\begin{equation}
\begin{cases}
 x'''(t)-\lambda [f(t,x(t)-\phi_\ld(t))+M]=0, &t\in(0,1),\\[.3pc]
 x(0)=x'(\eta)=x''(1)=0. &
\end{cases}
\end{equation}
Notice that
\begin{equation*}
x^*(t)\g \|x^*\|q(t)\g r'q(t)>\ld MKq(t)\g\phi_\ld(t),\quad \forall t\in (0,1).
\end{equation*}
Let $x(t)=x^*(t)-\phi_\ld(t)$ for $t\in [0,1]$, then $x(t)\in Q$. By
(16) or (18), it is easy to see that $x(t)$ is a positive solution to
$(1_\ld)$.

Similar to the above discussion, we can get another positive solution
$y^*$ to (19) from $\{y_j(t)\}_{j\g {j_0}}$ and $y^*(t)>\phi_\ld(t)$ for
$\fa t\in (0,1).$ Consequently, $y(t)=y^*(t)-\phi_\ld (t)$ is also a
solution to $(1_\ld)$.

We now show that (19) has no positive solution on $\partial Q_r$.
Suppose there exists $x\in \partial Q_r$ satisfying (19). Then
\begin{equation*}
x'''(t)\l \ld g(t)h(x(t)-\phi_\ld(t)), \quad \fa t\in (0,1).
\end{equation*}

Similar to the proof of Lemma 2.3, we can get $\ld \g\frac{b}{c(r)a}$.
This is in contradiction with $0<\ld < \frac{b}{c(r)a}$. Since
$r'\l\|x^*\|\l r$, $r\l\|y^*\|\l R$, then $x^{*}\neq y^*$. Consequently
$x\neq y$ and Theorem~2.1 follows.}\hfill $\Box$
\end{pot}

\begin{rema}
When $(1_{\ld})$ reduces to positone problem or $f(t,x)$ has no
singularity, Theorem~2.1 can also be used.
\end{rema}

\section{Examples}

\begin{exampl}
{\rm Consider the following BVP:
\begin{equation}
\begin{cases}
 x'''(t)-\ld f(t,x)=0, &t\in(0,1);\\[.3pc]
 x(0)=x'\left(\frac{2}{3}\right)=x''(1)=0, &
\end{cases}
\end{equation}
where
\begin{equation*}
f(t,x) = \frac{1}{\sqrt{t(1-t)}}\left(\frac{1}{x}+x^{2}\right)-\sin t.
\end{equation*}
It is easy to see $f(t,x)$ is singular at $t=0, t=1,$ and $x=0$, also
may be negative for some values of $t$ and $x$. Evidently, $(\hbox{H}_1)$ and
$(\hbox{H}_2)$ are satisfied for (20). Then by Theorem 2.1, we can get that (20) has
at least two positive solutions when $\ld$ is sufficiently small.}
\end{exampl}

\begin{exampl}
{\rm Consider the following BVP:
\begin{equation}
\begin{cases}
 x'''(t)-\ld f(t,x)=0, &t\in(0,1);\\[.3pc]
 x(0)=x'\left(\frac{3}{4}\right)=x''(1)=0, &
\end{cases}
\end{equation}
where
\begin{equation*}
f(t,x) = t\left(\frac{a}{x^{\a}} + b\hbox{e}^{x}\right),\quad a>0,\quad \a>0,\quad b>0.
\end{equation*}
It is easy to see that by Theorem~2.1, (21) has at least two positive
solutions when $\ld$ is sufficiently small.}
\end{exampl}

\section*{Acknowledgement}

This project was supported by the National Natural Science Foundation of
P.R.~China (10171057) and by the Natural Science Foundation of Shandong
Province (Z2000A02).

\end{document}